\documentclass[12pt,amsfonts]{article}
\usepackage{graphicx}
\usepackage{latexsym}
\usepackage{amssymb}
\usepackage{amsmath}
\usepackage{enumerate}
\usepackage{color}
\usepackage{layout}
\usepackage{dsfont}
\usepackage{eufrak}
\newtheorem{prop}{Proposition}
\newtheorem{lemma}{Lemma}

\newtheorem{theorem}{Theorem}
\newtheorem{remark}{Remark}

\def\real{{\mathord{{\rm I\kern-2.8pt R}}}}        
\def\inte{{\mathord{{\rm I\kern-2.8pt N}}}}

\def\sZZ{{\rm Z\kern-2.8ptem{}Z}}

\def\z{{\mathchoice
  {\sZZ}
  {\sZZ}
  {\rm Z\kern-0.30em{}Z}
  {\rm Z\kern-0.25em{}Z} }}
\def\sQQ{{\kern 0.27em \vrule height1.45ex width0.03em depth0em
          \kern-0.30em \rm Q}}
\def\qu{{\mathchoice
    {\sQQ}
    {\sQQ}
  {\kern 0.225em \vrule height1.05ex width0.025em depth0em \kern-0.25em \rm Q}
  {\kern 0.180em \vrule height0.78ex width0.020em depth0em \kern-0.20em \rm Q}
        }}
\def\sCC{{\kern 0.27em \vrule height1.45ex width0.03em depth0em
          \kern-0.30em \rm C}}
\def\complex{{\mathchoice
    {\sCC}
    {\sCC}
  {\kern 0.225em \vrule height1.05ex width0.025em depth0em \kern-0.25em \rm C}
  {\kern 0.180em \vrule height0.78ex width0.020em depth0em \kern-0.20em \rm C}
        }}

\newcommand{\R}{\mathbb{R}}

\newcommand{\ba}{\begin{array}}
\newcommand{\ea}{\end{array}}
\newcommand{\be}{\begin{equation}}
\newcommand{\ee}{\end{equation}}
\newcommand{\bea}{\begin{eqnarray}}
\newcommand{\eea}{\end{eqnarray}}
\newcommand{\beaa}{\begin{eqnarray*}}
\newcommand{\eeaa}{\end{eqnarray*}}

\newcommand{\eps}{\varepsilon}

\def\b{\beta}

\def\z{\zeta}

\font\tenmath=msbm10 \font\sevenmath=msbm7 \font\fivemath=msbm5
\newfam\mathfam \textfont\mathfam=\tenmath
\scriptfont\mathfam=\sevenmath \scriptscriptfont\mathfam=\fivemath
\def\math{\fam\mathfam}
\def \b{\noindent}

\def \={{\buildrel {\rm (law)} \over =}}

\def \R{{\math R}}

\def\qed{ \hfill \vrule width.25cm height.25cm depth0cm\smallskip}

\newcommand{\basa}{\begin{assumption}}
\newcommand{\easa}{\end{assumption}}

\newcommand{\bas}{\begin{assum}}
\newcommand{\eas}{\end{assum}}

\newcommand{\ignore}[1]{}
\begin{document}

\renewcommand{\thefootnote}{\fnsymbol{footnote}}

\renewcommand{\thefootnote}{\fnsymbol{footnote}}

\title{ Existence and Besov regularity of the density for a class of SDEs with Volterra noise  }
\author{Christian Olivera $^{1,}$\footnote{Supported by
FAPESP 2012/18739-0} \hskip0.2cm 
Ciprian A. Tudor $^{2}$ \vspace*{0.1in} \\
$^{1}$ Departamento de Matem\'atica, Universidade Estadual de Campinas,\\
13.081-970-Campinas-SP-Brazil. \\
colivera@ime.unicamp.br \vspace*{0.1in} \\
 $^{2}$ Laboratoire Paul Painlev\'e, Universit\'e de Lille 1\\
 F-59655 Villeneuve d'Ascq, France.\\
 \quad tudor@math.univ-lille1.fr\vspace*{0.1in}}
\maketitle

\begin{abstract}
By using a simple method based on the fractional integration by parts, we prove the existence and the Besov regularity of the density for solutions to stochastic differential equations driven by an additive  Gaussian Volterra process. We assume weak regularity  conditions on the drift.  Several examples of Gaussian Volterra noises are discussed.
\end{abstract}

\medskip


 \smallskip


\section{Introduction}

A new and simple method has been introduced in \cite{Deb2}, \cite{Fou} in order to obtain the absolute continuity of the law of  random variables.   In particular this method, based on fractional integration by parts, allows to obtain the existence of the density of solutions to stochastic differential equations (SDEs in the sequel), together with its Besov regularity, under low regularity assumptions on the coefficients of the equation. These new techniques avoid the use of the Malliavin calculus, which requires strong regularity of the coefficients of the SDE.  We refer, among others,  to \cite{Deb}, \cite{Deb2}, \cite{Fou}, \cite{Ro}, \cite{Sanz1}, \cite{Sanz2} for several applications of the fractional integration by parts methodology to concrete examples. 

Our purpose is to employ this new method in order to treat the case of SDE with additive Volterra noise, which has not yet been considered, as far as we know.  We consider the SDE in $\mathbb{R} ^{d}$
\begin{equation}
\label{1-intro}
X_{t}= x+ \int_{0} ^{t} b(s, X_{s}) ds+ B_{t}
\end{equation}
with $x\in \mathbb{R} ^{d}, b\in L^{\infty} ([0,T], C_{b}^{\beta } (\mathbb{R} ^{d} ) )$ and $(B_{t})_{t\in [0,T]}$ a $d$-dimensional Gaussian Volterra process that can be expressed as a Wiener integral with respect to the Wiener process under the form (\ref{volt}). Although our toy example is when $B$ is a $d$-dimensional fractional Brownian motion (fBm), we will show that many other examples of Volterra noises can be considered.

We will show that any strong solution to (\ref{1-intro}), when it exists, it admits a density with respect to the Lebesgue measure. Moreover,  we give the Besov regularity of the density of the solution, i.e. we find the Besov space to which the density belongs. Our main results are obtained under rather general condition on the noise (the class of examples includes the fractional Brownian motion and the  Ornstein-Uhlenbeck process, among others), and under a non-Lipschitz conditions on the drift $b$, i.e. $b\in L^{\infty} ([0,T], C_{b}^{\beta } (\mathbb{R} ^{d} ) )$. This method also works for the case of path-dependent SDEs.

We organized our paper as follows. In Section 2 we describe our context and our main assumptions. In Section 3 we prove the existence and the Besov regularity of the density of the solution to the stochastic differential equation (\ref{1-intro}). In Section 4 we extend our result to the path-dependent case. Section 5 contains several examples of Gaussian Volterra noises that fit our assumptions.

As a final remark on the notation: by $\vert \cdot \vert $ we denote the Euclidean norm in $\mathbb{R} ^{d}$, $C_{b}^{\alpha} $ denotes the set of bounded H\"older continuous functions of order  $\alpha$, while $C$ denotes throughout the paper a generic strictly positive constant that may change from line to line.

\section{Preliminaries}

We present below the basic definitions and assumptions.

\subsection{The context}

Let $(W_{t}) _{t\in [0, T]} = \left( W ^{(1)}_{t},..., W ^{(d)}_{t}\right) _{t\in [0, T]} $ be a $d$-dimensional Wiener process on the probability space $(\Omega, \mathcal{F}, P)$. Denote by $(\mathcal{F}_{t})_{t\in [0, T]}$ the filtration generated by $W$ and consider a Gaussian Volterra process $(B_{t}) _{t\in [0, T]} = \left( B ^{(1)}_{t},..., B^{(d)}_{t}\right) _{t\in [0, T]} $ that can be expresses as 
\begin{equation}
\label{volt}
B_{t}= \int_{0} ^ {t}K(t,s) dW_{s} 
\end{equation}
i.e. $B^{(i)}_{t}= \int_{0} ^ {t}K(t,s) dW^{(i) }_{s} $ for every $i=1,.., d$.  We assume in the sequel that $K$ is deterministic kernel such that 
\begin{equation*}
\int_{0} ^{T} K ^{2} (t,s) ds <\infty. 
\end{equation*}

We will the  consider the following SDE in $\mathbb{R} ^{d}$
\begin{equation}
\label{1}
dX_{t}= b(t, X_{t}) dt + B _{t} 
\end{equation}
with initial condition $x=(x_{1},.., x_{d})\in \mathbb{R} ^{d}$, where $(B_{t}) _{t\geq 0} $ is a Gaussian Volterra process of the form (\ref{volt}), i.e. for every $i=1,..,d$, 
\begin{equation*}
X_{t} ^{(i)}= x_{i}+ \int_{0}^{t} b_{i}(s, X_{s}) ds +B_{t} ^{(i)}
\end{equation*}
where $b_{i}$ are the components of the function $b$. We will assume that the drift coefficient in (\ref{1}) satisfies
\begin{equation}
\label{hyp1}
b\in L ^{\infty} \left( [0,T], C_{b} ^{\beta}(\mathbb{R}^{d}, \mathbb{R} ^{d} )\right) \mbox{ with } 0<\beta \leq 1.
\end{equation}

Notice that there is not a general result on the existence and uniqueness of the solution to (\ref{1}) under the assumption (\ref{hyp1}) for general Volterra noise of the form (\ref{volt}). In the sequel we will  work under the assumption that there exists a strong solution to (\ref{1}). Nevertheless, as we will comment in the last section,  there  are concrete situations when it exists an unique strong solution to (\ref{1})  under the assumptions (\ref{hyp1}) (for instance, this happens at least when the noise is a Wiener process or a fractional Brownian motion). 

If we assume stronger assumption on $b$ (i.e, that the drift is Lipschitz continuous and satisfies a linear growth condition, ), then we can easily get the existence and uniqueness of a strong solution to (\ref{1})  for a rather general Volterra noise $B$. In this case, the existence of the density of the solution to (\ref{1}) can be also obtained by different techniques (i.e. via Malliavin calculus). The advantage  of the method employed below is that it allows to find the Besov regularity of the density.

\subsection{Besov spaces}
We refer to \cite{Triebel} for a complete exposition on Besov spaces. Here we only recall the definition of a particular Besov space, namely the space $\mathcal{B} _{1, \infty}^{s} $ with $s>0$.

Consider a function $f:\mathbb{R} ^{d}\to \mathbb{R}$ and for every $x,h\in \mathbb{R}^{d}$, put
$$ (\Delta _{h} ^{1} f) (x)= f(x+h) -f(x)$$
 and for $n\geq 1$ integer, define the $n$th increment of the function $f$ at lag $h$

\[
(\Delta_{h}^{n}f )(x)= \Delta_{h}^{1}(\Delta_{h}^{n-1}f)(x)= \sum_{j=0}^{n} (-1)^{n-j}  f(x+jh). 
\]

\b For $0<s<n$ we define the norm
\begin{equation}\label{besovn}
\Vert f\Vert _{\mathcal{B} _{1,\infty}^{\alpha}}=\Vert f\Vert _{L ^{1}(\mathbb{R} ^{d})}+ \sup_{\vert h\vert \leq 1} \vert h\vert ^{-s} \Vert \Delta _{h} ^{n} f\Vert _{ L ^{1}(\mathbb{R} ^{d})}.
\end{equation}
It can be shown that for any $n,m>s$, the norms obtained in (\ref{besovn}) using $n,m$ are equivalent. Therefore, one can define the Besov space $\mathcal{B} _{1, \infty}^{s} $ as the set of functions $f\in L ^{1} (\mathbb{R} ^{d} )$ such that 
$$\Vert f\Vert _{\mathcal{B} _{1,\infty}^{\alpha}}<\infty.$$

\subsection{Fractional integration by parts}

Our main tool to get the existence and the regularity of the density of the solution to (\ref{1}) is the following smoothing lemma 
from \cite{Ro}.

\begin{lemma}\label{Ro}
 Let $X$ be a $\R^{d}$-valued random variable. If there exist
an integer $m\geq 1$, two  real numbers $s>0, \alpha>0$, with 
$\alpha <s <m$, and a constant $K>0$ such that for every $\phi\in C_{b}^{\alpha}(\R^{d})$ and 
$h\in \R^{d}$, with $|h|\leq 1$,

\[
\mathbf{E}[\Delta_{h}^{m}\phi(X)]\leq K |h|^{s} \|\phi\|_{C_{b}^{\alpha}}, 
\]
then $X$ has density $f_{X}$ with respect to Lebesgue measure on $\R^{d}$. Moreover, $f_{X}\in B_{1,\infty}^{s-\alpha}$
and

\[
\| f\|_{B_{1,\infty}^{s-\alpha}}\leq C (1+ K). 
\]

\end{lemma}

\section{The existence and the Besov regularity of the density}
We consider the setup from Section 2: the SDE (\ref{1}) with Volterra noise of the form (\ref{volt})  and with drift coefficient satisfying (\ref{hyp1}). We assume that there exists a strong solution to (\ref{1}). 

We fix a deterministic function  $\varphi \in C _{b} ^{\alpha} (\mathbb{R} ^{d} )$ wit  $\alpha \in (0,1)$  to be chosen later.   We need to estimate the quantity $\mathbf{E}  \left[ \Delta ^{m} _{h} \varphi (X _{t} ) \right]$ for $h>0$ and $m\geq 1$ integer. 

\smallskip
The core idea is to use   the auxiliary process

\begin{equation}\label{aux}
Y_{s}^{\epsilon} =\left \{
\begin{aligned}
    & X_{s},     \ \ \  \ \  \ \ \ \ \ \  \ \ \ \  \ \ s\leq t-\epsilon \,  \\
    \\[5pt]
		& X_{t-\epsilon} +   \int_{t-\epsilon}^{s} b(r,X_{t-\epsilon}) dr  + (B_{s}-B_{t-\epsilon}),  \ \ \  \ \ s\geq t-\epsilon \, .
		\end{aligned}
\right .
\end{equation}

\b We will  write
\begin{equation*}\mathbf{E}  \left[ \Delta ^{m} _{h} \varphi (X _{t} ) \right]= {\bf Pe}+ {\bf Ae}
\end{equation*}
where {\it  the probability estimate} {\bf Pe} is given by 

\begin{equation}
\label{pe}
{\bf Pe}= \mathbf{E} \left[ \Delta ^{m} _{h} \varphi (Y _{t} ^{\eps} ) \right]
\end{equation}
and {\it the approximation error } {\bf Ae} is  
\begin{equation}
\label{ae}
{\bf Ae}= \mathbf{E}  \left[ \Delta ^{m} _{h} \varphi (X _{t} ) \right]- \mathbf{E} \left[ \Delta ^{m} _{h} \varphi (Y _{t} ^{\eps} ) \right].
\end{equation}

We will deal  separately with the summands {\bf Pe} and {\bf Ae}, by using the ideas from \cite{Ro} and the properties of the Volterra noise $B$.

\subsection{The probabilistic estimate}
To get a suitable estimate for {\bf Pe}, we will express is in terms of two independent random variables. 
First notice that from (\ref{aux}), 
\begin{equation*}
Y _{t} ^{\eps} =X_{t-\eps} +\int_{t-\eps} ^{t} b(u, X_{t-\eps} ) du + B_{t}-B _{t-\eps}
\end{equation*}
and by writing $B_{t}- B_{t-\eps} = \int_{0} ^{t-\eps} \left( K(t,s)- K(t-\eps, s) \right) dW_{s}+  \int_{t-\eps}^{t} K(t,s) dW_{s}$, we obtain
\begin{eqnarray}
Y _{t} ^{\eps}
&=&X_{t-\eps} +\int_{t-\eps} ^{t} b(u, X_{t-\eps} ) du +\int_{0} ^{t-\eps} \left( K(t,s)- K(t-\eps, s) \right) dW_{s}+ \int_{t-\eps}^{t} K(t,s) dW_{s}\nonumber \\
&=& Z_{t}^{\eps} + I_{t} ^{\eps}\label{5m-1}
\end{eqnarray}
where
\begin{equation}
\label{ze}
Z_{t} ^{\eps}=X_{t-\eps} +\int_{t-\eps} ^{t} b(u, X_{t-\eps} ) du+ \int_{0} ^{t-\eps} \left( K(t,s)- K(t-\eps, s) \right) dW_{s}
\end{equation}
and $I_{t}^{i, \eps}=(I_{t} ^{1,\eps },..., I_{t} ^{d, \eps}) $ with 
\begin{equation}
\label{ie}
I_{t} ^{i, \eps}= \int_{t-\eps}^{t} K(t,s) dW_{i, s} \mbox{ for every } i=1,..,d.
\end{equation}
The key observation is that $Z_{t} ^{\eps}$  is  a $\mathcal{F}_{t-\eps}$ measurable  random variable in $\mathbb{R} ^{d}$  while $I_{t} ^{\eps} $  is a centered Gaussian random variable independent of$\mathcal{F}_{t-\eps}$.  Using the above decomposition (\ref{5m-1}), we obtain the following estimate for the probabilistic estimate

\begin{prop}\label{p1}
Assume (\ref{hyp1}) and suppose that for every real $h>0$ and for every integer $m\geq 1$
\begin{equation}
\label{cc1}
Var (I_{t}^{\eps})\geq C \eps ^{2A} K(\eps, t) \mbox{ with some } A \in (0,1)
\end{equation}
where $0<K(\eps, t)\leq C$ for every $\eps <t$. Then
$${\bf Pe} \leq C \Vert \varphi\Vert _{\infty}  \left( \frac{\vert h\vert }{\eps ^{A}}\right) ^{m} .$$
\end{prop}
{\bf Proof: } From the decomposition (\ref{5m-1}), with $\varphi  \in C _{b} ^{\alpha} (\mathbb{R} ^{d} )$,
\begin{eqnarray}
{\bf Pe}&=&\mathbf{E} \left[ \Delta ^{m} _{h} \varphi (Y _{t} ^{\eps} ) \right]= \mathbf{E} \left[ \Delta _{h} ^{m} \varphi (Z_{t} ^{\eps} +I_{t}^{\eps}) \right]\nonumber \\
&=& \mathbf{E} \left[ \mathbf{E}\left[ \Delta _{h} ^{m} \varphi (Z_{t} ^{\eps} +I_{t}^{\eps}) \right]/\mathcal{F} _{t-\eps} \right]= \mathbf{E} f( Z_{t} ^{\eps})\label{18m-1}
\end{eqnarray} 
with
\begin{equation*}
f(y)= \mathbf{E}\left[ \Delta _{h} ^{m} \varphi (y+ I_{t} ^{\eps} ) \right].
\end{equation*}
Denote by $g_{t, \eps} $ the density of the Gaussian random variable $I_{t}^{\eps} $, i.e.
\begin{equation}
\label{g}
g_{t, \eps} (x)= \frac{1}{\sqrt{2\pi Var (I_{t} ^{\eps}) } ^{d} } e ^{-\frac{ \vert x\vert ^{2}}{2Var(I_{t} ^{\eps} )}}.
\end{equation}
 We compute $f(y)$ via a trivial change of variables
\begin{eqnarray*}
f(y)&=&\int_{\mathbb{R} ^{d} } \Delta _{h} ^{m} \varphi (y+x) g_{t, \eps} (x)dx =\int_{\mathbb{R} ^{d}} \varphi (y+x) \left( \Delta _{-h} ^{m} g_{t, \eps }(x) \right) dx \\
&\leq & \Vert \varphi\Vert _{\infty} \Vert \Delta _{-h} ^{m} g_{t, \eps }(x) \Vert _{L ^{1} (\mathbb{R} ^{d} ) }.
\end{eqnarray*}
It follows from  \cite{Ro} that assumption (\ref{cc1}) implies that 
\begin{equation}\label{18m-2}
 \Vert \Delta _{-h} ^{m} g_{t, \eps (x) }\Vert _{L ^{1} (\mathbb{R} ^{d} ) }\leq  C \left(  \frac{\vert h\vert }{\eps ^{A}}\right) ^{m} 
\end{equation}
for every $h>0$ and for any integer $m\geq 1$. Then the conclusion is obtained from (\ref{18m-1}) and (\ref{18m-2}). \qed

We will see in the last section that (\ref{cc1}) is satisfied for many Gaussian processes, including the fractional Brownian motion.

\subsection{The approximation error}
In order to handle the term {\bf Ae} given by (\ref{ae}), we need the following hypothesis on the Gaussian noise $B$:  there exists $C>0$ such that

\begin{equation}
\label{cc2}
\mathbf{E} \left| B_{t}-B_{s} \right| ^{2} \leq C \vert t-s\vert ^{2H} \mbox{ with some } H\in (0,1).
\end{equation}

\begin{remark}
In particular, assumption (\ref{cc2}) implies that the process $B$ has H\"older continuous paths of order $\delta$ for every $\delta \in (0,H).$
\end{remark}

We have the following result for the approximation error {\bf Ae}.

\begin{prop}\label{p2}
Assume (\ref{hyp1}) and  (\ref{cc2}). Then  for every $0<\eps <t$,

\begin{equation}
\label{5m-3}
{\bf Ae}\leq C  \Vert \varphi \Vert _{ C_{b} ^{\alpha} } \eps ^{ (\beta H +1) \alpha}.
\end{equation}
\end{prop}
{\bf Proof: } Since $\varphi$ is $\alpha$-H\"older continuous, clearly

$$
{\bf Ae}= \mathbf{E}  \left[ \Delta ^{m} _{h} \varphi (X _{t} ) \right]- \mathbf{E} \left[ \Delta ^{m} _{h} \varphi (Y _{t} ^{\eps} ) \right] \leq \Vert \varphi \Vert _{ C_{b} ^{\alpha} }\mathbf{E} \left| X _{t}- Y_{t} ^{\eps} \right| ^{\alpha}. 
$$
Now, the difference  $X _{t}- Y_{t} ^{\eps}$ can be written as
\begin{equation*}
 X _{t}- Y_{t} ^{\eps}= \int_{t-\eps} ^{t} \left( b(u, X_{u})- b(u, X_{t-\eps} )\right) du. 
\end{equation*}
Thus
\begin{eqnarray*}
\mathbf{E} \left| X _{t}- Y_{t} ^{\eps} \right| ^{\alpha}&=& \mathbf{E} \left| \int_{t-\eps} ^{t} \left( b(u, X_{u})- b(u, X_{t-\eps} )\right) du\right| ^{\alpha} \nonumber\\
&\leq  &C \left| \int_{t-\eps} ^{t} \mathbf{E} \left |X_{u}-X_{t-\eps} \right| ^{\beta} du \right|  ^{\alpha}. 
\label{5m-2}
\end{eqnarray*}
Using (\ref{cc2}), for every $u>t-\eps$
$$\mathbf{E} \left |X_{u}-X_{t-\eps} \right| ^{\beta} =\mathbf{E} \left| \int_{ t-\eps} ^{u} b(v, X_{v}) dv+ B_{u}-B_{t-\eps} \right| ^{\beta } \leq  C \left( (u-t+\eps) ^{\beta}+ (u-t+\eps) ^{\beta H} \right).$$
So, by plugging the above inequality into (\ref{5m-2}),  
$$\mathbf{E} \left| X _{t}- Y_{t} ^{\eps} \right| ^{\alpha}\leq C \left| \int_{t-\eps} ^{t}  \left( (u-t+\eps) ^{\beta}+ (u-t+\eps) ^{\beta H} \right)du \right| ^{\alpha} \leq C \eps ^{ (\beta H +1) \alpha}$$
and this implies (\ref{5m-3}). \qed

\subsection{The density of the solution}

We are now ready to apply the smoothing Lemma \ref{Ro}. From Proposition \ref{p1} and \ref{p2} we obtain:

\begin{theorem}\label{t1}
Assume (\ref{hyp1}), (\ref{cc1}) and (\ref{cc2}). Let $(X_{t}) _{t\in [0,T]} $ be a strong solution to (\ref{1}). Then for every $t\in [0, T]$, the random variable $X_{t}$ admits a density $\rho_{t}$ with respect to the Lebesgue measure. Moreover,
$$\rho_{t} \in \mathcal{B} _{1, \infty} ^{\eta}  \mbox{ for any } \eta < \frac{1-A+\beta H}{A}.$$
\end{theorem}
{\bf Proof: } From Propositions \ref{p1} and \ref{p2}
$$\mathbf{E}\Delta _{h} ^{m} \varphi(X_{t}) \leq {\bf Pe}+ {\bf Ae} \leq C\Vert \varphi \Vert _{ C_{b} ^{\alpha} }   \left(\left(  \frac{\vert h\vert }{\eps ^{A}}\right) ^{m}+ \eps ^{ (\beta H +1) \alpha}\right).$$
Let us choose
$$\eps= h ^{\frac{m}{\alpha (\beta H+1)+Am}}.$$
Then we get 
\begin{equation*}
\mathbf{E}\Delta _{h} ^{m} \varphi(X_{t}) \leq C \Vert \varphi \Vert _{ C_{b} ^{\alpha} }  \vert h\vert ^{s}  
\end{equation*}
with
$$s= \frac{m\alpha(1+\beta H)}{\alpha (1+\beta H)+Am}.$$

Note that for $m$ large enough, the exponent of $\vert h\vert $ is about $\frac{\alpha (1+\beta H)}{A}$. Therefore, by Lemma \ref{Ro},  for every $t\in (0, T]$, the random variable $X_{t}$ has a density $\rho_{t}$ belonging to the Besov space $\mathcal{B} _{1,\infty} ^{\eta} $, with $\eta< s-\alpha= \frac{\alpha (1+\beta H)}{A}-\alpha$. Since we can choose $\alpha $ to be arbitrary close to 1, we obtain the conclusion. \qed

Let us finish this section but some comments around Theorem \ref{t1}.

\begin{remark}\label{rem1}
\begin{itemize}

\item In the case of the Wiener noise (i.e. $K(t,s)= 1_{[0, t]}(s)$ for every $s,t\in [0, T]$, conditions (\ref{cc1}) and (\ref{cc2}) hold with $A=H=\frac{1}{2}$.  On the other hand, a unique strong solution to (\ref{1}) exists under (\ref{hyp1}). Indeed,  the existence and uniqueness of the strong solution is assured  for every measurable function $b\in L ^{\infty} ([0,T]\times \mathbb{R} ^{d})$ (see \cite{Zvo} for $d=1$ and \cite{Ve} for general dimensional $d\geq 1$). It follows from Theorem \ref{t1}, that the solution to (\ref{1}) admits a density in the Besov space $\mathcal{B} _{1, \infty}^{\eta}$ for every $\eta <1+\beta.$ We retrieve a result in Section 2 of \cite{Ro}.

\item
We notice that both the noise in (\ref{1}) and the variance of $I_{t} ^{\eps}$ affect the regularity of the density. More regular are the paths of the noise $B$ (i.e. $H$ increases), more regular is the density of solution (i.e. $\eta $ incresases). Also, as the variance of $I_{t} ^{\eps} $ increases, then $A$ decreases and therefore the regularity of the solution increases.

\end{itemize}
\end{remark}

\section{The path dependent case}

The argument from the previous section can be easily adapted to treat the path dependent case. As before, we will consider $(W_{t}) _{t\in [0, T]} $ a $d$-dimensional $\mathcal{F}_{t}$- Brownian motion on the probability space $(\Omega, \mathcal{F}, P)$ and let $(B_{t})_{t\in [0,T] }$ be a Volterra process of the form (\ref{volt}). We consider the SDE

\begin{equation}
\label{sde2}
X_{t} (x)= x+ \int_{0}^ {t} b(r,V_{r}, X_{r} ) dr +  B_{t}
\end{equation}
with $t\in [0,T]$, $x\in \mathbb{R} ^ {d}$. In this section, the drift $b$ is assumed to satisfy

\begin{equation}\label{H2}
b\in L^{\infty}([0,T], C_{b}^{\beta} (\R^{d} \times \R^{d},\R^{d} ), \ \ \ 0<\beta\leq 1,
\end{equation}
while $(V_{t}) _{t\in [0,T]} $ is a $\mathcal{F}_{t}$-adapted process such that 
\begin{equation}\label{H3}
\mathbf{E} |V_{t}-V_{s}|^{\beta}\leq C |t-s|^{\delta} \mbox{ for some } \delta >0.
\end{equation}

We assume, as before, that there exists a strong solution to (\ref{sde2}). For $\eps>0$, we define the auxiliary process $Y_{t} ^{\eps}$ by

\begin{equation}\label{aux2}
Y_{s}^{\epsilon} =\left \{
\begin{aligned}
    & X_{s},     \ \ \  \ \  \ \ \ \ \ \  \ \ \ \  \ \ s\leq t-\epsilon \,  \\
    \\[5pt]
		& X_{t-\epsilon} +   \int_{t-\epsilon}^{s} b(r,V_{t-\epsilon},X_{t-\epsilon}) dr  + (B_{s}-B_{t-\epsilon}),  \ \ \  \ \ s\geq t-\epsilon \, .
		\end{aligned}
\right .
\end{equation}

We decompose again the quantity $\mathbf{E} [\Delta_{h}^{m}\varphi(X_{t})]$ into two terms, the approximation error  

\begin{equation}\label{error}
{\bf Ae}=\mathbf{E}[\Delta_{h}^{m}\varphi(X_{t})]-\mathbf{E}[\Delta_{h}^{m}\varphi(Y_{t}^{\epsilon})]
\end{equation}
and the probabilistic estimated  

\begin{equation}\label{prob}
{\bf Pe }=\mathbf{E} [\Delta_{h}^{m}\varphi(Y_{t}^{\epsilon})].
\end{equation}

\b Concerning the summand {\bf Pe}, we have the following estimate:

\begin{lemma}\label{Pe3} Assume (\ref{H2}) and  (\ref{cc1}). Then we have 

\begin{equation}
\label{18m-4}
{\bf Pe  }=\mathbf{E} [\Delta_{h}^{m}\varphi(Y_{t}^{\epsilon})]\leq  C  \| \varphi  \|_{\infty}  \left(\frac{\vert h\vert }{\eps ^{A}}\right) ^{m}.
\end{equation}
\end{lemma}
{\bf Proof: } From (\ref{aux2}), we can write
$${\bf Pe}= Z_{t}^{\eps} + I_{t}^{\eps}$$
where $I_{t}^{\eps}$ is given by (\ref{ie}) and 
\begin{equation*}
Z_{t}^{\epsilon}=X_{t-\epsilon} +  \int_{t-\epsilon}^{t} b(r,V_{t-\epsilon},X_{t-\epsilon})  ds   + \int_{0}^{t-\epsilon} \left( K(t,s)-K(t-\epsilon,s)\right)  dB_{s}. 
\end{equation*}
since $Z_{t} ^{\eps}$ is $\mathcal{F}_{t-\eps} $ measurable and $I_{t} ^{\eps}$ is independent by $\mathcal{F}_{t-\eps} $, we can write

\begin{eqnarray*}
\mathbf{E} [\Delta_{h}^{m}\varphi(Y_{t}^{\epsilon})]&=& \mathbf{E} [\Delta_{h}^{m}\varphi(Z_{t}^{\epsilon}+ I_{t} ^{\eps} )]\\
&=& \mathbf{E} \big[\mathbf{E} [\Delta_{h}^{m}\varphi(y+I_{t} ^{\eps} )]_{y=Z^{\epsilon}}\big].
\end{eqnarray*}
and using (recall that $g_{t, \eps}$ is given by (\ref{g})) 

\begin{eqnarray*}
\mathbf{E} [\Delta_{h}^{m}\varphi(y+I_{t} ^{\eps} )]&=& \int_{\R^{d}} \varphi(y+x) \Delta_{-h}^{m}g_{t, \epsilon}(x)  dx 
\\
&
\leq  &\| \varphi \|_{\infty}  \|\Delta_{-h}^{m}g_{t, \epsilon}(x)  \|_{L^{1} (\mathbb{R} ^{d}} 
\leq  C_{H} \| \varphi \|_{\infty}  (\frac{\vert h\vert }{\epsilon^{A}})^{m}.
\end{eqnarray*}
\qed

For the approximation error term {\bf Ae}, we have the next result.

\begin{prop}
Assume (\ref{H2}), (\ref{H3}) and (\ref{cc2}).  Then 

\begin{equation}\label{Ae3}
{\bf Ae  }\leq C_{H} \|\varphi \|_{\alpha} \epsilon^{ \alpha (\mu+1)}
\end{equation}
where $\mu=min(\beta H, \delta)$.

\end{prop}
{\bf Proof: } We write as in the proof of Proposition \ref{p2}
\begin{eqnarray}
Ae&=&\mathbf{E} [\Delta_{h}^{m}\varphi(X_{t})]-\mathbf{E} [\Delta_{h}^{m}\varphi(Y_{t}^{\epsilon})]\nonumber \\
&\leq&
\Vert \varphi \Vert _{C_{b} ^{\alpha}} \mathbf{E} |\int_{t-\epsilon}^{t} (b(r,V_{r},X_{r})-b(r,V_{t-\epsilon}, X_{t-\epsilon}))  dr   |^{\alpha}\nonumber \\
&\leq & 
\| \varphi \|_{\alpha} \mathbf{E} |\int_{t-\epsilon}^{t} |b(r,V_{r},X_{r})-b(r,V_{r}, X_{t-\epsilon})| + 
  |b(r,V_{r},X_{t-\epsilon})-b(r,V_{t-\epsilon}, X_{t-\epsilon})| dr|^{\alpha}
\nonumber \\
&
\leq &  \|\varphi \|_{\alpha}  \  \big( \mathbf{E} \int_{t-\epsilon}^{t} |X_{r}-X_{t-\epsilon}|^{\beta}+  |V_{r}-V_{t-\epsilon}|^{\beta}   dr \big)^{\alpha}. \label{5m-4}
\end{eqnarray}

By insert the following two bounds 
\[
 \mathbf{E} |X_{r}-X_{t-\epsilon}|^{\beta} \leq  \| b \|_{L^{\infty}} (r-t+\epsilon)^{\beta}+ (r-t+\epsilon)^{\beta H}
\]

and 

\[
 \mathbf{E} |V_{r}-V_{t-\epsilon}|^{\beta} \leq  C (r-t+\epsilon)^{\delta}
\]
into (\ref{5m-4}),  we get 

\begin{equation*}
Ae\leq \|\varphi \|_{\alpha} \big(\int_{t-\epsilon}^{t} (r-t+\epsilon)^{\beta}+ (r-t+\epsilon)^{\beta H} + (r-t+\epsilon)^{\delta}  dr \big)^{ \alpha }
\leq C_{H} \|\varphi \|_{\alpha} \epsilon^{ \alpha (\mu+1)}
\end{equation*}
where $\mu=min(\beta H, \delta)$. \qed

\begin{theorem} \label{t2}We assume the conditions (\ref{cc1}), (\ref{cc2}), (\ref{H2}) and (\ref{H3}). Then 
 the law of $X_{t}$ has density $\rho_{t,x}$ respect to the Lebesgue measure and 
$\rho\in B_{1,\infty}^{\eta}$  with $\eta< \frac{\mu+1-A}{A}$ where $\mu=min(\delta,H\beta)$. 
\end{theorem}
{\bf Proof: } From  the estimates (\ref{Ae3}) and (\ref{18m-4}) we get

\[
\mathbf{E} [\Delta_{h}^{m}(X_{t})]\leq  C_{H} \|\varphi \|_{\alpha} \epsilon^{\alpha(\mu+1)} +  C_{H} \| \varphi \|_{\alpha}  (\frac{\vert h\vert }{\epsilon^{A}})^{m} 
\]

Now, choosing $\epsilon=h^{\frac{m}{\alpha (\mu +1)+Am}}$ and proceeding as in the proof of Theorem \ref{t1}, we obtain the desired conclusion. \qed

\begin{remark}
Notice that the Besov regularity of the density is affected by the regularity of the process $V$ since the exponent $\delta$  from (\ref{H3}) appears in the above result. By taking a regular process $V$ with $\delta >H\beta$, we retrieve the result in  Theorem \ref{t1}, but for a process $V$ such that $\delta <H\beta$, the Besov regularity of the density will change. 

\end{remark}

\section{Examples}
We discuss several examples where our main results stated in Theorems \ref{t1} and \ref{t2} applies.

\subsection{Fractional Brownian motion}

Let $(B_{t}) _{t\in [0, T]}$ be a fractional Brownian motion with Hurst index $H\in (0,1)$. Recall that $B$ is a centered Gaussian process with covariance
\begin{equation*}
  \mathbf{E} B_{t} B_{s}= \frac{1}{2} \left( t ^{2H}+ s^{2H}-\vert t-s\vert ^{2H} \right) \mbox{ for every }s,t\in [0,T].
\end{equation*}
The fBm admits the following integral representation
\begin{equation}
\label{fbm1}
B_{t}= \int_{0} ^{t} K_{H}(t,s) dW_{s}
\end{equation}
where $(W_{t})_{t\in [0, T]}$ is a Wiener process, and $K_{ H}(t,s)$ is the kernel
\begin{equation}
K_{H}( t,s)=d_{H}\left(  t-s\right)  ^{H-\frac{1}{2}}+s^{H-\frac
{1}{2}}F_{1}\left(  \frac{t}{s}\right)  , \label{for1}%
\end{equation}
$d_{H}$ being a constant and
\begin{equation*}
F_{1}\left(  z\right)  =d_{H}\left(  \frac{1}{2}-H\right)  \int _{0}^{z-1}\theta^{H-\frac{3}{2}}\left( 1-\left(
\theta+1\right)
^{H-\frac{1}{2}}\right)  d\theta.
\end{equation*}
If $H>\frac{1}{2}$, the kernel $K_{H}$ has the simpler expression
\begin{equation}
 \label{K}
K_{H}(t,s)= c_{H} s^{\frac{1}{2}-H} \int _{s}^{t} (u-s)^{H-\frac{3}{2}} u^{H-\frac{1}{2}}  du
 \end{equation}
 where $t>s$ and $c_{H} =\left( \frac{ H(H-1) }{\beta( 2-2H, H-\frac{1}{2}) } \right)
^{\frac{1}{2}}.$

The SDE (\ref{1}) with fBm noise has been treated in \cite{NO1}, \cite{NO2}, \cite{Ngu}, among others. The following facts have been proven  (for $d=1$):

\begin{itemize}
\item If $H>\frac{1}{2}$, then there exists a unique strong solution to (\ref{1}) if the drift $b$ is H\"older continuous in time of order $\gamma > H-\frac{1}{2}$ and it is H\"older continuous in space of order $\alpha > 1-\frac{1}{2H}$, i.e. 
$$\vert b(t,x)-b(s,y)\vert \leq C (\vert x-y\vert ^{\alpha} + \vert t-s\vert ^{\gamma} )$$
with $\alpha > 1-\frac{1}{2H}$ and $\gamma > H-\frac{1}{2}$.

\item If $H<\frac{1}{2}$, then there exists a unique strong solution to (\ref{1}) if $b$ satisfies the linear growth condition 
\begin{equation}
\label{cl}
\vert b(t,x) \vert \leq C (1+\vert x\vert)
\end{equation}
for every $t\in [0, T], x\in \mathbb{R}$.

\item  If $H=\frac{1}{2}$, see Remark \ref{rem1}.

\end{itemize}

Notice that the assumption (\ref{hyp1}) clearly implies the linear growth condition (\ref{cl}), thus we always have existence and uniqueness of the solution to (\ref{1}) under (\ref{hyp1}). When $H>\frac{1}{2}$, we need to assume $\beta \geq 1-\frac{1}{2H} $ in (\ref{hyp1}) and also
$$\vert b(t,x)-b(s,x) \vert\leq C \vert t-s\vert ^{\gamma} $$
with $\gamma > H-\frac{1}{2}$ for every $s,t\in [0, T], x\in \mathbb{R}$.

In order to apply Theorems \ref{t1} and \ref{t2}, we need to check (\ref{cc1}) and (\ref{cc2}). Assumption (\ref{cc2}) clearly holds with $H=\frac{1}{2}$. To check (\ref{cc1}), we discuss separately the cases $H>\frac{1}{2}$ and $H<\frac{1}{2}$. 

\vskip0.2cm 

\b {\bf  The case  $H>\frac{1}{2}$. } We  claim  that   $ I_{t}^ {\epsilon}=\int_{t-\epsilon}^{t} K_{H} (t,s) dB_{s}$ is 
    Gaussian  with expectation equal to zero and variance  bigger than $  c_{H} \epsilon^{2H}$.  We have 
\[
Var(I_{t}^{\epsilon})= \mathbf{E}  | \int_{t-\epsilon}^{t} K_{H} (t,s) dB_{s}   |^{2} =\int_{t-\epsilon}^{t} |K_{H}(t,s)|^{2}  ds 
\]
and from formula (\ref{K}), since for $H-\frac{1}{2}>0$ we have $\left(\frac{u}{s}\right) ^{H -\frac{1}{2}}\geq 1$, we can write

\[
K_{H}(t,s)\geq c_{H}  \int _{s}^{t} (u-s)^{H-\frac{3}{2}}   du=C_{H} (t-s)^{H-\frac{1}{2}}.
\]
Then

\[
\int_{t-\epsilon}^{t} |K(t,s)|^{2}  ds\geq C_{H}^{2} \int_{t-\epsilon}^{t}  (t-s)^{2H-1} ds =C_{H} \epsilon^{2H}
\]
so (\ref{cc1}) holds with $A=H$.

\vskip0.3cm

\b {\bf The case $H<\frac{1}{2}$. } From Proposition 5.12 in \cite{N} we have 

\[
K(t,s)\geq C_{H} (\frac{t}{s})^{H-\frac{1}{2}} \ (t-s)^{H-\frac{1}{2}}.
\]

\b Thus

\begin{eqnarray*}
\int_{t-\epsilon}^{t} |K(t,s)|^{2}  ds &\geq& C_{H}  t^{2H-1}   \int_{t-\epsilon}^{t}  s^{1-2H} (t-s)^{2H-1} ds\\
&\geq &
C_{H} t^{2H-1} (t-\epsilon)^{1-2H} \int_{t-\epsilon}^{t}  (t-s)^{2H-1} ds
\\
&=&
C_{H} t^{2H-1} (t-\epsilon)^{1-2H}  \epsilon^{2H}\\
&=&
C_{H} t^{2H-1} (1-\frac{\epsilon}{t})^{1-2H}  \epsilon^{2H}.
\end{eqnarray*}
Consequently (\ref{cc1}) holds with $A=H$ and $K(\eps, t)=  t^{2H-1} (1-\frac{\epsilon}{t})^{1-2H} $ which is less that a constant for $0<\eps <t$.

\subsection{The Riemann-Liouville process}

The Riemann-Liouville process is defined as 
\begin{equation}
\label{rl}
B_{t}=\int_{0} ^{t} (t-s) ^{H-\frac{1}{2}} dW_{s}, \mbox{ for every } t\in [0, T]
\end{equation}
with $H\in (0,1)$. It shares many properties with the fBm (it is self-similr of index $H$, is paths are H\"older continuous of order $\delta \in (0,H)$), but it has not stationary increments. Notice that 
$$Var (I_{t} ^{\eps} ) = \int_{t-\eps} ^{t}(t-s) ^{2H-1} ds= \frac{1}{2H} \eps ^{2H}$$ 
and it is well-known that 
$$\mathbf{E} \left| B_{t}- B_{s} \right| ^{2}\leq C\vert t-s\vert ^{2H}.$$
Therefore assumptions (\ref{cc1}) and (\ref{cc2}) are fulfilled with $A=H$ and $K(t,\eps)=1$.

\subsection{The Ornstein-Uhlenbeck process}
The Ornstein-Uhlenbeck process $(B_{t}) _{t\in [0,T]}$ can be expressed as 
$$B_{t}= \int_{0} ^{t} e ^{-(t-s) }dW_{s}.$$
It represents the unique solution to the SDE $dB_{t}=-B_{t}dt +dW_{t}$ with vanishing initial condition.  It is well-known that 
$$\mathbf{E}  \left| B_{t}- B_{s} \right| ^{2}\leq C\vert t-s\vert $$
so (\ref{cc2}) is satisfied with $H=\frac{1}{2}$. On the other hand, if $I_{t} ^{\eps}$ is given by (\ref{ie}),
$$Var (I_{t} ^{\eps} ) =Var \left( \int_{t-\eps} ^{t}  e ^{-(t-s) }dW_{s}\right) = \int_{t-\eps} ^{t} e ^{-2(t-s) }ds =\frac{1}{2} (1-e ^{-2\eps} )\geq c\eps$$
so (\ref{cc1}) holds with $A=\frac{1}{2}$ and $K(\eps, t)=1$.

\end{document}